\documentclass[12pt]{article}
\usepackage{amsmath,amssymb,amsfonts,amsthm}
\newenvironment{myabstract}{\par\noindent
{\bf Abstract . } \small }
{\par\vskip8pt minus3pt\rm}
\newcounter{item}[section]
\newcounter{kirshr}
\newcounter{kirsha}
\newcounter{kirshb}
\newenvironment{enumroman}{\setcounter{kirshr}{1}
\begin{list}{(\roman{kirshr})}{\usecounter{kirshr}} }{\end{list}}
\newenvironment{enumarab}{\setcounter{kirshb}{1}
\begin{list}{(\arabic{kirshb})}{\usecounter{kirshb}} }{\end{list}}
\newtheorem{theorem}{Theorem}[section]

\newtheorem{lemma}[theorem]{Lemma}

\theoremstyle{definition}

\newtheorem{definition}[theorem]{Definition}

\def\C{{\mathfrak{C}}}

\def\Nr{{\mathfrak{Nr}}}

\def\Sg{{\mathfrak{Sg}}}

\def\A{{\mathfrak{A}}}
\def\B{{\mathfrak{B}}}
\def\C{{\mathfrak{C}}}

\def\M{{\mathfrak{M}}}

\def\CA{{\bf CA}}

\def\RCA{{\bf RCA}}

\def\(R)RA{{\bf (R)RA}}

 \def\CA{{\sf CA}}
\def\B{{\sf B}}

\def\Nr{{\mathfrak{Nr}}}

\def\Nr{{\mathfrak{Nr}}}

\def\A{{\mathfrak{A}}}
\def\B{{\mathfrak{B}}}
\def\C{{\mathfrak{C}}}

\def\A{{\mathfrak{A}}}
\def\B{{\mathfrak{B}}}
\def\C{{\mathfrak{C}}}

\def\L{{\mathfrak{L}}}

\def\L{{\mathfrak{L}}}
\def\Bl{{\mathfrak{Bl}}}
\def\CA{{\bf CA}}

\def\RCA{{\bf RCA}}

\def\L{{\mathfrak L}}

\def\Nr{{\mathfrak{Nr}}}

\def\CA{{\bf CA}}
\def\RCA{{\bf RCA}}

\def\Sg{{\mathfrak Sg}}
\def\P{{\mathfrak P}}

\title{Which fragments of $\L_{\infty}$ is the class of neat reducts sensitive to?}
\author{Tarek Sayed Ahmed \\
Department of Mathematics, Faculty of Science,\\ 
Cairo University, Giza, Egypt.
  }
\begin{document}
\maketitle

\begin{myabstract} Let $\L$ be a quantifier predicate logic. Let $K$ be a class of algebras.
We say that $K$ is sensitive to $\L$, if there is an algebra in $K$, that is $\L$ interpretable into an another algebra, 
and this latter algebra is elementary equivalent to an algebra not
in $K$. (In particular, if $\L$ is $\L_{\omega,\omega}$, this means that $K$ is not elementary).

Let $L$ be the full typless logic studied in [HMT2] with $\alpha$ many variables. Let $\L$ be the quantifier free reduct 
of first order logic, endowed with infinite conjunctions and possibly constants.
We shall construct a relativized model $M$ in $L$ and a cylindric set algebra $\A_M$ of dimension 
$\alpha$, such that $\A_M\in \Nr_{\alpha}\CA_{\alpha+\omega}$ and algebras $\B, \C$ such that $\A$ is $\L$ interpretable in $\B$,
$\C$ is elementary equivalent to $\B$
and $B\notin \Nr_{\alpha}\CA_{\alpha+1}$. In particular, the class of neat reducts is sensitive to $\L$. 

\end{myabstract}

\section{Introduction}

The class of neat reducts has been studied extensively by the author and others. 
This class is not closed under forming subalgebras. Several conditions strengthening forming ordinary subalgebras, taking {\it strong} subalgebras, 
in a certain precise, sense have been suggested, 
so that the resulting algebra is a neat reduct.  The question, can be paraphrased as follows: Can we define a part of the neat reduct, that is of
course closed under 
the operations, in a certain logic that forces this part to be also a neat reduct. When we ask about ordinary subalgebras, then we are in equational logic.
For example is a complete subalgebra of a neat reduct a neat reduct. 
The answer is no. But here we have a flavour of second order logic, so the question is 
how sensitive is the operation of forming subalgebras to other logics. 

It is possible, in theory that for example elementary subagebras are neat reducts.
Or perhaps subalgebras that satisfy the same $L_{\kappa, \omega}$ sentences, or $L_{\infty}$ sentences, are neat reducts
A lot of investigations have shown that the class of neat reducts  is truly reselient to being closed under various 
kinds of substructures, or strong subalgebras. This has also proved to be surprising, in another sense, for such classes resultng 
from taking special sub neat reducts have turned out closely related to other natural notions, appearing
in completely different contexts, like for example the class of completely representable algebras, 
which in turn is closely related to the metalogical property of omiting types
and also to  various amalgamation bases for representable algebras. The latter is closely related to the interpolation property. 
These metalogical properties are viewed in fragments
of the full typless logic.

So there is a spectrum of kinds of subalgebras, and such a spectrum is strongly related to various 
metalogical properties of the corresponding logic.

In this paper, we show that the class of neat reducts is also sensitive to a natural fragment of $L_{\infty}$, namely 
quanifier free reduct of first order logic  endowed with infinite conjunctions.

We alert the reader that we will deal with three logics. Two are based in advance. Standard first order logic and  
more basic predicate logic which is the full logic 
studied in  [HMT] sec 4.3, and finally, a transient logic $\L$, which is actually a parameter, in all cases it is a  reduct of $L_{\infty}$.
In all logics we have an infinite supply of variables and 
we have equality interpreted in the intended models the usual way.

Full logics are based on relational languages, the relations have a fixed arity, specified by an infinite ordinal $\alpha$.
Atomic formulas are of the form $R(x_0, x_1,\ldots )$ so that the variables can only occur in their natural order.
Such formulas are called restricted.  We have the usual Boolean connectives, namely $\land$ and $\neg.$ 
Quantification is only allowed on finitely many variables, and formulas are specified recursively the usual way.
On the other hand, $\L$ could be $L_{\omega,\omega}$ itself, or $L_{\kappa.\omega}$, 
$\kappa$ a regular cardinal, or $L_{\infty}$, or a reduct thereof.

Let $\A$ be a first order structure and $\B$ be an $\L$ structure. Then $\A$ is $\L$ interpretable in $\B$,
(we wil be concerned only with one dimensional quantifier free interpretations) 
if there is a function $f:\A\to \B$, such that for any formula $\phi$ of first order logic, 
there exists an $\L$ formula $\psi$ such that for all $a\in \A$, 
we have $$\A\models \phi(a)\Leftrightarrow  \B\models\psi(f(a)).$$ 
When $\L$ is first order logic this is the usual definition of interpretability, which is a generalization of the notion of relativized reducts.

For full logics, we will consider relativized semantics. While full logics are relational, we only allow function symbols in $\L$, we do not have
relation symbols, so that atomic formulas are of the form $t=s$, where $t$ and$s$ are terms.
While models for full logic will be relational, for $\L$ we consider only 
algebras with possibly infinitary operations.


We will construct a model for the full language, that generalizes Fraisse's constructions, 
of ordinary models, particulary those which have elimination of 
quantifiers. Our proof is a non trivial step-by- step construction, that can be implemented using games.
The meta logic used here in our construction is ordinary first order logic.

Then we will define a weak cylindric set algebra on this model; this algebra is in $\Nr_{\alpha}\CA_{\alpha+\omega}$. 

The question is: For which such transiant logics $\L$, and $k\in \omega$
is $\A$,  $\L$ interpretable in an algebra that is elementary equivalent to an algebra not in $\Nr_{\alpha}\CA_{\alpha+k}$.
The answer is known when $\L$ is usual first order logic,  here we investigate the analogous situation for the quantifier free reduct of $L_{\infty}$.

\begin{definition} 
Let $\alpha$ be an ordinal 
and $M$ be a set. A weak space of dimension $\alpha$ and base $M$
is a set of the form
$$\{s\in {}^{\alpha}M: |\{i\in \alpha: s_i\neq p_i\}| <\omega\}$$
for some $p\in {}^{\alpha}M.$ 
We denote this set by $^{\alpha}M^{(p)}.$
Let $\Lambda$ be a full language 
having $\beta$ many relation symbols..
A weak structure for $\Lambda$ is a triple ${\M}= (M, R, p)$ where
$M$ is a non empty set $p\in {}^{\alpha}M$ and 
$R$ is a function with domain $\beta$ 
assigning to each $i<\beta$ a subset $R_i^{\M}$ of the weak space
$^{\rho i}M^{(p)}.$ 
A sequence $s\in {}^{\alpha}M^{(p)}$ 
satisfies an atomic formula 
$R_i(v_0, \ldots v_i\ldots)_{i<\rho(i)}$
if $s \in R_i^{\M}$. $s$ satisfies $v_i=v_j$ if $s(i)=s(j).$ 
We can extend the notion of 
satisfiability to all formulas in the usual Tarskian way. 
\end{definition}
$R$ be an uncountable set and let $Cof R$ be set of all non-empty finite or cofinite subsets  $R$.
Let $\alpha$ be an ordinal. For $k$ finite, $k\geq 1$, let
$$S(\alpha,k)=\{i\in {}^\alpha(\alpha+k)^{(Id)}: \alpha+k-1\in Rgi\},$$
$$\eta(X)=\bigvee \{C_r: r\in X\},$$
$$\eta(R\sim X)=\bigwedge\{\neg C_r: r\in X\}.$$


Let $(W_i: i\in \alpha)$ be a disjoint family of sets each of cardinality $|R|$.
Let $M$ be their disjoint union, that is
$M=\bigcup W_i$. Let $\sim$ be an equvalence relation on $M$ such that $a\sim b$ iff $a,b$ are in the same block.
Let $T=\prod W_i$. Let $s\in T$, and let $V={}^{\alpha}M^{(s)}$. 
For $s\in V$, we write $D(s)$ if $s_i\in W_i$, and we let $\C=\wp(V)$.
 
\begin{lemma}
There are $\alpha$-ary relations $C_r\subseteq {}^{\alpha}M^{(s)}$ on the base $M$ for all $r\in R$,
such that conditions (i)-(v) below hold:
\begin{enumroman}
\item $\forall s(s\in C_r\implies D(s))$

\item For all $f\in {}^{\alpha}W^{(s)}$ for all $r\in R$, for all permutations
$\pi\in ^{\alpha}\alpha^{(Id)}$, if $f\in C_r$ then $f\circ \pi\in C_r.$ 

\item For all $1\leq k<\omega$, for all 
$v\in {}^{\alpha+k-1}W^{(s)}$ one to one,  for all $x\in W$, $x\in W_m$ say, then for any
function $g:S(\alpha,k)\to Cof^+R$ 
for which $\{i\in S(\alpha,k):|\{g(i)\neq R\}|<\omega\}$, 
there is a $v_{\alpha+k-1}\in W_m\smallsetminus Rgv$ such that 
and  
$$\bigwedge \{D(v_{i_j})_{j<\alpha}\implies \eta(g(i))[\langle v_{i_j}\rangle]: 
i\in S(\alpha,k)\}.$$
\item The $C_r$'s are pairwise disjoint.
\end{enumroman}
\end{lemma}

For $u\in S_{\alpha}$ and $r\in R$, let 
$$p(u,r)=C_r\cap (W_{u_0}\times W_{u_1}\times W_{u_i}\times)\cap {}^{\alpha}W^{(s)}.$$ 
Let $$\A=\Sg^{\C}\{p(u,r): u\in S_{\alpha}: r\in R\}.$$
Then $\A$ is the weak set algebra based on $\M$. 

For $u\in {}^{\alpha}\alpha^{(Id)}$, let $1_u=W_{u_0}\times W_{u_i}\times \cap V$
and $\A_u$ denote the relativisation of $\A$ to $1_u$.
i.e $$\A_u=\{x\in A: x\leq 1_u\}.$$ 

$\A_u$ is a boolean algebra. Also  $\A_u$ is uncountable and atomic for every $u\in V$
The sets $C_r^{\M}$, for $r\in R$ are disjoint 
elements of $\A_u$.  
Because of the saturation condition above, we have $\A\in \Nr_{\alpha}\CA_{\alpha+\omega}$.

Define a map $f: \Bl\A\to \prod_{u\in {}V}\A_u$, by
$$f(a)=\langle a\cdot \chi_u\rangle_{u\in{}V}.$$

Now consider the langauge $\L$.
We will expand the language of the boolean algebra $\prod_{u\in V}\A_u$ by constants in 
such a way that
$\A$ becomes interpretable in the expanded structure.
For this we need. Let $\P$ denote the 
following structure for the signature of boolean algebras expanded
by constant symbols $1_u$ for $u\in {}V$ and ${\sf d}_{ij}$ for $i,j\in \alpha$: 

\begin{enumarab}
\item The boolean part of $\P$ is the boolean algebra $\prod_{u\in {}V}\A_u$,

\item $1_u^{\P}=f(\chi_u^{\M})=\langle 0,\cdots0,1,0,\cdots\rangle$ 
(with the $1$ in the $u^{th}$ place)
for each $u\in {}V$,

\item ${\sf d}_{ij}^{\P}=f({\sf d}_{ij}^{\A})$ for $i,j<\alpha$.
\end{enumarab}

Define a map $f: \Bl\A\to \prod_{u\in {}V}\A_u$, by
$$f(a)=\langle a\cdot \chi_u\rangle_{u\in{}V}.$$

Here, and elsewhere, for a relation algebra $\C$, $\Bl\C$ denotes its boolean reduct.

We now show that $\A$ is $\L$ interpretable in $\P.$  
For this it is enough to show that 
$f$ is one to one and that $Rng(f)$ 
(Range of $f$) and the $f$-images of the graphs of the cylindric algebra functions in $\A$ 
are definable in $\P$. Since the $\chi_u^{\M}$ partition 
the unit of $\A$,  each $a\in A$ has a unique expression in the form
$\sum_{u\in {}V}(a\cdot \chi_u^{\M}),$ and it follows that 
$f$ is boolean isomorphism: $bool(\A)\to \prod_{u\in {}V}\A_u.$
So the $f$-images of the graphs of the boolean functions on
$\A$ are trivially definable. 
$f$ is bijective so $Rng(f)$ is 
definable, by $x=x$. For the diagonals, $f({\sf d}_{ij}^{\A})$ is definable by $x={\sf d}_{ij}$.

Finally we consider cylindrifications for $i<\alpha$. Let $S\subseteq {}V$ and  $i,j<\alpha$, 
let $t_S$ and $h_S$ be the closed infinitary terms:
$$\sum\{1_v: v\in {}V, v\equiv_i u\text { for some } u\in S\}.$$
 
Let
$$\eta_i(x,y)=\bigwedge_{S\subseteq {}V}(\bigwedge_{u\in S} x.1_u\neq 0\land 
\bigwedge_{u\in {}V\smallsetminus S}x.1_u=0\longrightarrow y=t_S).$$
These are well defined.

We claim that for all $a\in A$, $b\in P$, we have 
$$\P\models \eta_i(f(a),b)\text { iff } b=f({\sf c}_i^{\A}a).$$
To see this, let $f(a)=\langle a_u\rangle_{u\in {}V}$, say. 
So in $\A$ we have $a=\sum_ua_u.$
Let $u$ be given; $a_u$ has the form $(\chi_i\land \phi)^{\M}$ for some $\phi\in L^3$, so
${\sf c}_i^A(a_u)=(\exists x_i(\chi_u\land \phi))^{\M}.$ 
By property (vi), if  $a_u\neq 0$, this is 
$(\exists x_i\chi_u)^M$; by property $5$,
this is $(\bigvee_{v\in {}V, v\equiv_iu}\chi_v)^{\M}.$
Let $S=\{u\in {}V: a_u\neq 0\}.$
By normality and additivity of cylindrifications we have,
$${\sf c}_i^A(a)=\sum_{u\in {}V} {\sf c}_i^Aa_u=
\sum_{u\in S}{\sf c}_i^Aa_u=\sum_{u\in S}(\sum_{v\in {}V, v\equiv_i u}\chi_v^{\M})$$
$$=\sum\{\chi_v^{\M}: v\in {}V, v\equiv_i u\text { for some } u\in S\}.$$
So $\P\models f({\sf c}_i^{\A}a)=t_S$. Hence $\P\models \eta_i(f(a),f({\sf c}_i^{\A}a)).$
Conversely, if $\P\models \eta_i(f(a),b)$, we require $b=f({\sf c}_ia)$. 
Now $S$ is the unique subset of $V$ such that 
$$\P\models \bigwedge_{u\in S}f(a)\cdot 1_u\neq 0\land \bigwedge_{u\in {}V\smallsetminus S}
f(a)\cdot 1_u=0.$$  So we obtain 
$$b=t_S=f({\sf c}_i^Aa).$$
Choose any countable boolean elementary 
subalgebra of $bool(A_{Id})$, $B_{Id}$ say.
Thus $B_{Id}\preceq bool(A_{Id})$.
By the Feferman-Vaught theorem   
$$Q=((B_{Id}\times \prod_{u\in V\smallsetminus {Id}} 
bool(A_u)),1_u,d_{ij})_{u\in {}V, i,j<\alpha}\preceq$$  
$$(\prod_{u\in {}V} bool(A_u)),1_u, d_{ij})_{u\in V}=P.$$
Let $\B$ be the result of applying the interpretation given above to $Q$.
Then $\B\equiv \A$ as cylindric algebras. 
Therefore $\B\in \RCA_{\alpha}$. Finally, $\B\notin \Nr_{\alpha}\CA_{\alpha+1}.$ 
\end{document}